\documentclass[11pt,a4paper,reqno]{amsart} 

\usepackage{amssymb,graphicx}

\usepackage{amssymb}

\newcommand{\koniec}{\begin{flushright}  $\Box $ \end{flushright}}
\newtheorem{theo}{Theorem}[section]
\newtheorem{prop}[theo]{Proposition}

\newcounter{mnotecount}[section]

\renewcommand{\themnotecount}{\thesection.\arabic{mnotecount}}

\newcommand{\mnote}[1]
{\protect{\stepcounter{mnotecount}}$^{\mbox{\footnotesize
$
\bullet$\themnotecount}}$ \marginpar{
\raggedright\tiny\em
$\!\!\!\!\!\!\,\bullet$\themnotecount: #1} }

\newcommand{\CP}{\mathbb{CP}}
\newcommand{\C}{\mathbb{C}}
\newcommand{\spp}{\mathbb{S}}

\newcommand{\PP}{\mathbb{P}}

\newcommand{\R}{\mathbb{R}}

\def\p{\partial}
\def\OO{\mathcal{O}}

\def\ov{\overline}

\def\be{\begin{equation}}
\def\ee{\end{equation}}

\def\bea{\begin{eqnarray}}
\def\eea{\end{eqnarray}}

\def\ov{\overline}

\title{Jumps, folds,  and singularities of Kodaira moduli spaces.}

\author{Maciej Dunajski, James Gundry and Paul Tod}
\address{Maciej Dunajski. Department of Applied Mathematics and Theoretical Physics\\
University of Cambridge\\ Wilberforce Road, Cambridge CB3 0WA\\ UK.\\
Paul Tod. The Mathematical Institute\\
University of Oxford\\ Andrew Wiles Building \\
Woodstock Road, Oxford OX2 6GG, UK.}

\begin{document}\date{25 April 2018}
\maketitle
\begin{abstract}
For any integer $k$ we construct an explicit example of a twistor space
which contains a one--parameter family of jumping rational curves,
where the normal bundle changes from $\OO(1)\oplus\OO(1)$ to
$\OO(k)\oplus\OO(2-k)$. For $k>3$ the resulting anti--self--dual
Ricci--flat manifold is a Zariski cone in the space of holomorphic sections of
$\OO(k)$. In the case $k=2$ we recover the canonical example
of Hitchin's folded hyper--K\"ahler manifold, where the jumping lines form a
three--parameter family. We show that in this case there exist normalisable solutions to the
Schr\"odinger equation which extend  through the fold.

\end{abstract}
\section{Introduction}
The non--linear graviton twistor construction of
Penrose \cite{penrose}
gives a one-to-one correspondence between
holomorphic anti--self--dual Ricci flat metrics on complex four--manifolds $M_\C$, and complex three--folds $\mathcal{Z}$ with a family of
rational curves. The points in $M_\C$ correspond to
holomorphic sections of $\mathcal{Z}\longrightarrow \CP^1$ characterised by their normal bundle $\OO(1)\oplus\OO(1)$, where $\OO(k)$ is a line bundle
over $\CP^1$ with Chern class $k$.

If the normal bundles of rational curves corresponding to points on
a surface $S\subset M_\C$
(of co-dimension one or more) change, then the metric becomes singular on
$S$. There are other geometric structures, most notably the self--dual
two--forms spanning $\Lambda^{2}_+$, which nevertheless remain regular on
$S$. In the case of a single jump to $\OO(2)\oplus\OO$
this results in a folded hyper--K\"ahler structure in the sense of Hitchin
\cite{H2}. Examples of such structures,  and the underlying existence
theorems are known \cite{biqard}, and some applications in theoretical physics have recently emerged \cite{reall}.

The aim of this paper is to construct an explicit example of a twistor
space where the normal bundle jumps from $\OO(1)\oplus\OO(1)$ to
$\OO(k)\oplus\OO(2-k)$ for any integer $k> 2$. This jump occurs
on a curve $\gamma\subset M_\C$, and the corresponding metric
on $M_\C$ can be constructed explicitly. It admits a tri--holomorphic Killing
vector, and so away from the jump it can be put in the standard
Gibbons--Hawking form \cite{GH} by a coordinate transformation.
This transformation removes the region of $M_\C$ where the jump occurs.
The resulting metric on $M_\C\setminus\gamma$ is still singular
on the surface $S$ where the Gibbons--Hawking harmonic function vanishes.
The normal bundles of twistor lines corresponding to the points on $S$ jump to
$\OO(2)\oplus\OO$.

 In the next Section we shall set up the twistor correspondence, where
the non--deformed twistor space $\mathcal{Z}$ is an affine line bundle over the total space of $\OO(k)$, for any $k$. In Theorem \ref{main_theo}
we find  Kodaira deformations preserving
this affine bundle, and leading to a four--manifold $M_\C$ arising as
a Zariski cone in the $(k+1)$--dimensional space of holomorphic sections
$H^0(\CP^1, \OO(k))$.
In Sections \ref{k=3ex}  and \ref{k=4ex} we give expressions for the metric in cases
where $k=3$ and $k=4$, and show that the Gibbons-Hawking potential
\[
V(X, Y, Z)=\frac{1}{2(k-1)!}
\frac{\p^{k-1}}{\p \lambda^{k-1}}
(\lambda^2 (X-Y)+2\lambda Z+(X+Y))^{-1/2}|_{\lambda=0},
\]
on $\C^3$ or $\R^{1,2}$ corresponds to the general $k$.
In Section \ref{EHST} we give an explicit coordinate transformation 
between a real form of the Sparling--Tod solution and 
the Eguchi--Hanson gravitational instanton and its  dyonic limit.

In Section \ref{sec_casca} it will be shown how a cascade of intermediate jumps
\[
\OO(1)\oplus\OO(1)\rightarrow
\OO(2)\oplus\OO\rightarrow\dots\rightarrow\OO(k)\oplus\OO(2-k)
\]
arises on surfaces on $M_\C$ with various co-dimensions.
In Section \ref{sec_leg} we shall put the construction in the framework
of the generalised Legendre transform
\cite{IR96,B00,DM}, and show how
the Zariski cone $M_\C$ arises as the zero--locus of a zero-rest-mass
field corresponding to a cohomology class in $H^1(\OO(k), \OO(2-k))$.

Finally in Section \ref{sec_shrod} we shall come back
to the case $k=2$, where the metric arising form Theorem \ref{main_theo}
admits a Riemannian real slice $(M, g)$ , which
is the canonical model of a folded hyper--K\"ahler structure
\be
\label{folded_gg}
g=Z(dX^2+dY^2+dZ^2)+Z^{-1}\Big(dT +\frac{1}{2}XdY-
\frac{1}{2}YdX\Big)^2,
\ee
where $S\subset M$ given by $Z=0$ is the fold. Answering a question of Manton, we shall show that
despite the blow-up in the metric there exist normalisable solutions to the
Schr\"odinger equation which extend through the fold.
\subsection*{Acknowledgements.} We are grateful to Nigel Hitchin for discussions about
folded geometry, and to Nick Manton for his suggestion that normalisable solutions
to the Schr\"odinger equations may extend through the folds. The work of MD has been partially
supported by STFC consolidated grant no. ST/P000681/1. JG was supported by an STFC studentship. PT acknowledges the hospitality of Girton College, Cambridge, where he held the Brenda Ryman Visiting Fellowship while this work was underway, and CMS for office space and computing facilities.
\section{Twistor spaces as affine bundles}
We shall start off by reviewing the twistor 
correspondence \cite{penrose,HT,MD}.
Let $M_\C$ be a complex four--manifold with a holomorphic orientation vol, and a holomorphic
Ricci--flat metric $g$
such that the Weyl tensor is anti--self--dual (ASD). The anti--self--duality
is the Frobenius integrability condition for the existence of a three--parameter
family of self--dual totally null surfaces ($\alpha$--surfaces) in $M_\C$, and the twistor space
$\mathcal{Z}$ is the three--dimensional complex manifold with the $\alpha$--surfaces as points.
This leads to a double fibration picture
\[
M_\C\longleftarrow {\mathcal F} \stackrel{\rho}\longrightarrow \mathcal{Z},
\]
where ${\mathcal F}\subset M_\C\times \mathcal{Z}$ is the space of incident pairs $(p, \xi)$ such that
$p\in M_\C$ lies on an $\alpha$--surface $\xi\subset M_\C$.
 A point in $M_\C$ corresponds to a projective line $L_p\cong\CP^1$ in $\mathcal{Z}$ which consists of all
$\alpha$--surfaces through $p$. A conformal structure $[g]$ on $M_\C$ is encoded in the algebraic
geometry of curves in $\mathcal{Z}$: two points in $M_\C$ are null--separated iff the corresponding curves
in $\mathcal{Z}$  intersect in one point.

There are two additional structures on $\mathcal{Z}$ resulting from the existence of a Ricci--flat
metric $g\in [g]$, and the canonical isomorphism
\[
TM_\C=\spp\otimes\spp',
\]
where $\spp$ and $\spp'$ are two rank--two complex symplectic vector bundles over $M_\C$.
\begin{itemize}
\item
The Levi--Civita connection of $g$ gives a  flat spin connection\footnote{To avoid the repeated usage of primed spinor indices in Section \ref{sec_leg} we depart from the usual twistor conventions, and swap the roles of primed and unprimed indices.} on $\spp$. Thus there exists
a two--dimensional space of parallel sections of
$\spp$. This, together with the isomorphism
$\Lambda^2_+\cong \spp\odot\spp$ and a natural identification
${\mathcal F}=\PP(\spp)$, gives a holomorphic projection
\be
\label{projection_pt}
\mu:\mathcal{Z}\longrightarrow \CP^1,
\ee
such that the points in $M_\C$ are holomorphic sections of $\mu$ with normal bundle
$\OO(1)\oplus\OO(1)$.
\item The  parallel basis $(\Sigma^{00}, \Sigma^{01}, \Sigma^{11})$ of $\Lambda^2_+$
gives rise to a symplectic two--form $\Sigma$ on the fibres of
(\ref{projection_pt}) which takes values in the line bundle $\OO(2)$. The pull back of
$\Sigma$ from $\mathcal{Z}$ to $\PP(\spp)$ is
\[
\rho^*(\Sigma)=\Sigma^{11}-2\lambda \Sigma^{01}+\lambda^{2}\Sigma^{00},
\]
where $\lambda$ is an affine coordinate on the fibres of $\PP(\spp)$.
\end{itemize}
The two--form $\Sigma$ fixes a volume form  vol on $M_\C$: the condition $\Sigma\wedge\Sigma=0$
gives
\[
\mbox{vol}=\Sigma^{00}\wedge\Sigma^{11}=-2\Sigma^{01}\wedge\Sigma^{01}.
\]

\subsection{Jumping lines}
\label{sect_jump_lines}
Let $[\pi_{0}, \pi_{1}]$ be homogeneous coordinates on
$\CP^1$. Cover $\CP^1$ with two open sets
\be
\label{UUopen_sets}
U=\{[\pi]\in\CP^1, \pi_{1}\neq 0\}, \quad
\widetilde{U}=\{[\pi]\in\CP^1, \pi_{0}\neq 0\}
\ee
and set $\lambda=\pi_{0}/\pi_{1}$ on $U\cap\widetilde{U}$.
The Birkhoff-Grothendieck theorem  states that any rank-$2$ holomorphic vector bundle
over $\CP^1$ is isomorphic to a direct sum of line bundles
$\OO(p)\oplus\OO(q)
$
for some integers $p, q$. Moreover the transition matrix
$F:\C^*\rightarrow GL(2, \C)$ of this bundle can be written
as
\be
\label{BH_split}
F=\widetilde{H}\; \mbox{diag}(\lambda^{-p}, \lambda^{-q})\; H^{-1},
\ee
where $H:U\rightarrow GL(2, \C)$ and
$\widetilde{H}:\widetilde{U}\rightarrow GL(2, \C)$ are holomorphic.

Let $\mathcal{Z}_b\longrightarrow \CP^1$ be
a one-parameter family of
rank-two vector bundles determined by the patching matrix
\[
F_b=\left (
\begin{array}{cc}
\lambda^{k-2} & b{\lambda}^{-1} \\
0 & {\lambda^{-k}}
\end{array}
\right ),
\]
where $b$ is a constant and $k$ is a positive integer. If $b=0$
then $F_0=\mbox{diag}(\lambda^{k-2}, \lambda^{-k})$
is the patching matrix for $\mathcal{Z}_{0}=\OO(2-k)\oplus\OO(k)$ with $H$ and $\widetilde{H}$ in (\ref{BH_split}) both equal to the identity matrix.  If $b\neq 0$ then
\[
F_b=\left(
\begin{array}{cc}
0 & {b} \\
-b^{-1} & {\lambda}^{1-k}
\end{array}
\right)
\left(
\begin{array}{cc}
{\lambda}^{-1} & 0 \\
0 & {\lambda}^{-1}
\end{array}
\right )
\left(
\begin{array}{cc}
1 & 0 \\
-b^{-1}\lambda^{k-1} & 1
\end{array}
\right )^{-1}
\]
which is of the form (\ref{BH_split}).
Thus $\mathcal{Z}_b=\OO(1)\oplus\OO(1)$ if $b\neq 0$.
This is the twistor space with the holomorphic sections of $\mathcal{Z}_b\rightarrow \CP^1$ parametrised by points in $M_\C=\C^4$ with the flat metric\footnote{If $k=2$, and $b$ is interpreted as the inverse
of the speed of light, then the jumping from $\mathcal{Z}_b$ to $\mathcal{Z}_0$ is the Newtonian limit
of the twistor correspondence \cite{DG16}.}.
\subsection{Twistor space as an affine bundle over
$\OO{(k)}$} Let $(Q, \lambda)$ and $(\widetilde{Q},
\tilde{\lambda}=\lambda^{-1})$ be coordinates on the  pre-images of $U$ and $\widetilde{U}$ in the  total space
of the line bundle $\OO(k)\rightarrow \CP^1$.
On the pre-image of $U\cap\widetilde{U}$ in $\mathcal{Z}_b$ we have
\be
\label{patch_flat}
\widetilde{Q}={\lambda^{-k}}Q,\quad
\tilde{\tau}=\lambda^{k-2}\tau+b\lambda^{-1} Q,
\ee
where $\tau$ and $\tilde{\tau}$ are coordinates on the fibres of ${\mathcal Z}_b$ over $U$ 
and $\widetilde{U}$ respectively.
Restricting the inhomogeneous coordinates to a
section of $\OO(k)\rightarrow \CP^1$
\be
\label{q1}
Q=x_k\lambda^k+\dots +x_1\lambda+x_0
\ee
and performing the splitting of (\ref{patch_flat}) gives
\be
\label{split_2}
\tilde{\tau}-b\lambda^{-1}x_0-bx_1=\lambda^{k-2}\tau
+b\lambda^{k-1}x_k+\dots +b\lambda x_2.
\ee
The LHS and RHS of (\ref{split_2}) are holomorphic on $\widetilde{U}$ and $U$ respectively.
Both sides of this relation are sections of a line bundle with a negative Chern class, so should vanish by the Liouville theorem.
Thus
\be
\label{tau_1}
\tau=-b(\lambda x_k-x_{k-1}-\lambda^{-1}x_{k-2}-\dots-
\lambda^{k-3}x_2).
\ee
This will be holomorphic in $\lambda$ on $U$ if $(k-3)$ conditions
\be
\label{flat_conditions}
x_{k-2}=x_{k-3}=\dots=x_2=0
\ee
hold. These conditions arise only if $k>3$.
They define a holomorphic four--dimensional subspace
$M_\C$ in the $(k+1)$ dimensional space of holomorphic sections of
$\OO(k)$.

 The algebraic geometry of holomorphic sections of
$\mathcal{Z}_b\rightarrow\CP^1$ determines a conformal
structure on $M_\C$: two points in $M_\C$ are null separated iff the corresponding sections intersect at one point
in $\mathcal{Z}_b$. Infinitesimally, a vector in $T_pM_\C$ is null
if the corresponding section of $N(L_p)$ vanishes at one point. This condition is equivalent to the existence of the unique solution
$\lambda=\lambda_0$ to a simultaneous system
\be
\label{penrose_metric}
\delta Q=0, \quad \delta \tau=0,
\ee
where $Q$ and $\tau$ are given by (\ref{q1}) and
(\ref{tau_1}).
These conditions give
\[
\lambda^k\delta x_k+\dots +\lambda\delta x_1+\delta x_0=0, \quad
\lambda \delta x_k+\delta x_{k-1}=0.
\]
Imposing (\ref{flat_conditions})  and using the second equation replaces the first equation by
$\lambda\delta x_1+\delta x_0=0$. Eliminating $\lambda$
between the two equations in (\ref{penrose_metric}) gives
the quadratic conformal structure
\be
\label{final_flat}
[g]=\delta x_k\delta x_0-\delta x_{k-1}\delta x_1
\ee
which is flat.
\section{Conformal structures on Zariski cones
from Kodaira deformations.}
\label{section_zari}
Consider an affine line bundle $\mathcal{Z}\rightarrow \OO(k)$, with underlying translation bundle
given by $\OO(2 - k)$. Such bundles are classified by elements
of $H^1(\OO(k), \OO(2-k))$ and we choose a cohomology representative which leads to
the patching relations
\be
\label{patch_2}
\widetilde{Q}={\lambda^{-k}}Q,\quad
\tilde{\tau}=\lambda^{k-2}\tau+a\lambda^{-2}Q^2, \quad\mbox{where}\;
a=\mbox{const.}
\ee
Restricting this to  holomorphic sections (\ref{q1}) of $\mathcal{Z}\rightarrow\CP^1$
and splitting gives
\be
\label{split_one}
\tilde{\tau}-a\lambda^{-2}({x_0}^2+2\lambda x_0x_1+
\lambda^2(2x_0x_2+{x_1}^2))=
\lambda^{k-2}\tau+a\lambda(2x_0x_3+2x_1x_2)+\dots
+a\lambda^{2k-2} {x_k}^2.
\ee
Therefore $\tau$ is holomorphic in $\lambda$ if $(k-3)$ quadratic conditions
\be
\label{cones}
x_0x_3+x_1x_2=0, \quad x_0x_4+x_1x_3+\frac{{x_2}^2}{2}=0, \quad \dots, \quad
x_0x_{k-1}+x_1x_{k-2}+\dots\;\;\;
\ee
hold. These constraints put no restrictions on $x_k$, and we can assume that
\[
t=x_0, \quad z=x_1, \quad y=x_2, \quad x=x_k
\]
are coordinates on an open set in $M_\C$, and that the remaining coordinates $(x_3, \dots,  x_{k-1})$
have been expressed as functions of $(y, z, t)$. To compute the
ASD
conformal structure $[g]$ on $M_\C$
we follow the procedure leading to
(\ref{final_flat}), except that to  simplify the computations the condition
(\ref{penrose_metric})  is replaced by the equivalent condition ${\delta \widetilde{Q}}=0,
\delta \tilde{\tau}=0$ (the resulting conformal structure does not depend on the choice of the open set) and pull the differentials $\delta x_i$ in
${\delta \widetilde{Q}}$  back to $M_\C$. To eliminate
$\lambda$ we take the resultant Res of the quadratic $\lambda^2\delta\tilde{\tau}$ and the polynomial of degree $k$ given by
${\lambda^k\delta\widetilde{Q}}$. The resultant is a section of $\mbox{Sym}^{(k+2)}(T^*M_\C)$
which factorises as $[g] (\delta x_0)^k$, where $[g]\in \mbox{Sym}^{2}(T^*M_\C)$ is the conformal
structure given by
\be
\label{resultant_g}
[g]=\mbox{Res}\Big({x_0}\delta x_0+\lambda( x_0\delta x_1+ x_1\delta x_0) +
\lambda^2(x_0\delta x_2+x_2\delta x_0+x_1\delta x_1), \delta x_0 +\lambda \delta x_{1}+\dots+\lambda^k\delta x_k
\Big).
\ee
Here $\delta x_3, \dots, \delta x_{k-1}$  are the pullbacks from
$H^{0}(\CP^1, \OO(k))$ to $M_\C$ defined by the relations
(\ref{cones}). The deformation (\ref{patch_2}) preserves the fibration $\mu:\mathcal{Z}\rightarrow \CP^1$, and the fibres of $\mu$ are equipped with an $\OO(2)$--valued symplectic form
\[
\Sigma=\lambda^{-2}dQ\wedge d\tau=d\widetilde{Q}\wedge d\tilde{\tau}.
\]
Therefore there exists a Ricci--flat metric $g\subset[g]$ in the ASD conformal class
(\ref{resultant_g}).
\begin{theo}
\label{main_theo}
Let $(M_\C, g)$ be an ASD Ricci--flat manifold corresponding to the twistor space
with the patching relations (\ref{patch_2}).
There exists a curve $\gamma\subset M_\C$ such that all points on $\gamma$ correspond to rational curves  $\mathcal{Z}$ where the normal bundle jumps from $\OO(1)\oplus\OO(1)$
to $\OO(k)\oplus\OO(2-k)$, and such that $\gamma$ is preserved by
a tri--holomorphic Killing vector.
\end{theo}
{\bf Proof}.
We shall first prove that $(M_\C, g)$ admits a triholomorphic Killing vector field.
The coefficients of the conformal structure (\ref{resultant_g})
do not depend on $x\equiv x_k$, so $K=\p/\p x$ is
a conformal Killing vector. Conformal Killing vectors in $M_\C$ generate one-parameter groups of transformations
of $M_\C$ which map $\alpha$--surfaces to $\alpha$--surfaces. Thus (as the points in
$\mathcal{Z}$ are $\alpha$--surfaces in $M_\C$) conformal Killing vectors correspond to
global holomorphic vector fields on $\mathcal{Z}$. Consider the holomorphic vector field
$\mathcal{K}$
in $\mathcal{Z}$ corresponding to the conformal Killing vector  $K=\p/\p x$. We shall
compute this vector field on the open set $\widetilde{U}$
\begin{eqnarray*}
{\mathcal K}&=&\frac{\p \widetilde{Q}}{\p x}
\frac{\p}{\p \widetilde{Q}}+
\frac{\p \tilde{\tau}}{\p x}\frac{\p}{\p \tilde{\tau}}\\
&=&\frac{\p}{\p\widetilde{Q}}
\end{eqnarray*}
where we have used (\ref{split_one}). Therefore
\[
{\mathcal L}_{\mathcal K} \Sigma=0, \qquad {\mathcal L}_{\mathcal K} \lambda=0,
\]
and so ${\mathcal K}$ preserves  the symplectic two-form on the fibres of $\mathcal{Z}\rightarrow \CP^1$, as well as the fibration itself.
The first condition implies that $K$ is a Killing vector of the Ricci--flat metric
singled out by $\Sigma$ in the conformal structure $[g]$. The second condition
means that  $K$ is tri--holomorphic (it acts trivially on the basis of parallel self-dual two--forms).

Now consider the  normal bundle $N(L_p)$ to a curve $L_p\subset \mathcal{Z}$ corresponding to a point $p\in M_\C$.
For a generic $p$, the bundle $N(L_p)$ is biholomorphic to $\OO(1)\oplus\OO(1)$.
Its patching matrix is given by
\[
F_N=
\left(
\begin{array}{cc}
\frac{\p \tilde{\tau}}{\p\tilde \tau} &
\frac{\p \tilde{\tau}}{\p\tilde Q}\\
\frac{\p \widetilde{Q}}{\p\tilde \tau} &
\frac{\p \widetilde{Q}}{\p\tilde Q}
\end{array}
\right)=
\left(
\begin{array}{cc}
\lambda^{k-2} & 2a\lambda^{-2} Q  \\
0& \lambda^{-k}
\end{array}
\right ).
\]
To  investigate the non-generic points introduce the splitting matrices
\[
H=\left(
\begin{array}{cc}
1 & h  \\
0& 1
\end{array}
\right ),
\quad
\widetilde{H}=\left(
\begin{array}{cc}
1 & \tilde{h}\\
0& 1
\end{array}
\right )
\]
which are invertible and holomorphic in $U$ and $\widetilde{U}$ respectively.
Now
\[
\widetilde{H}F_N H^{-1}=\left(
\begin{array}{cc}
\lambda^{k-2} &
2a\lambda^{-2} Q
-h\lambda^{k-2}+\tilde{h}\lambda^{-k} \\
0& \lambda^{-k}
\end{array}
\right ).
\]
Thus the normal bundle jumps from $\OO(1)\oplus\OO(1)$ to $\OO(2-k)\oplus\OO(k)$
at points of $M_\C$ where
\[
2a\lambda^{-2} Q-h\lambda^{k-2}+\tilde{h}\lambda^{-k}=0.
\]
The functions $h$ and $\tilde{h}$ are holomorphic in $\lambda$ and $\lambda^{-1}$ respectively.
Therefore the equality can be satisfied by some choice of $h$ and $\tilde{h}$ only if
\[
x_{k-1}=x_{k-2}=\dots=x_0=0.
\]
These conditions imply the conditions (\ref{cones}), and leave the coordinate $x_k$ unspecified.
Thus there is a one--parameter family $\gamma\subset M_\C$ of jumping lines in $\mathcal{Z}$,
and $K$ is tangent to $\gamma$.
\koniec
In addition to $K$, the Ricci--flat metric arising from Theorem  \ref{main_theo} admits
a homothetic conformal
Killing vector $x\p_x+y\p_y+z\p_z+t\p_t$.
In what follows we shall work out the metrics and their Gibbons-Hawking forms in detail when
$k=3$ and $k=4$. We note that in these two cases the metric admits another Killing vector, where
the coordinates scale with different weights.
\subsection{Ricci--flat metric with $k=3$}
\label{k=3ex}
Parametrise the sections of $\OO(3)\rightarrow\CP^1$ by
\[
Q=t+\lambda z+\lambda^2 y+\lambda^3 x.
\]
In this case the anti-self-dual conformal structure (\ref{resultant_g}) is
\be
\label{g_2}
g=\Omega^2\Big(t^3 dx^2+2t^2z dx dy+t(yt+z^2)dxdz+
z(3yt-z^2)dxdt+tz^2 dy^2+z(yt+z^2)dydz
\ee
\[
+y(yt+z^2)dydt+yz^2 dz^2+
2y^2 z dzdt+y^3 dt^2\Big).
\]
We find that the choice
\be
\label{conf_f}
\Omega^2=\frac{2a^3}{z^2-yt}
\ee
makes the resulting metric Ricci--flat.
If the coordinates $(x, y, z, t)$ are chosen to be real, then
$g$ is real, and has neutral signature. The basis of self-dual
parallel  2-forms is
\begin{eqnarray*}
\Sigma^{00}&=&2dx\wedge (tdy+ydt+zdz),\\
\Sigma^{01}&=&tdx\wedge dz+zdx\wedge dt+zdy\wedge dz+ydy\wedge dt, \\
\Sigma^{11}&=&2(tdx+ydz+zdy)\wedge dt.
\end{eqnarray*}
These forms are Lie derived by the  Killing vector
$K=\p/\p x$. This Killing vector is tri-holomorphic
and therefore the metric $g$ can be cast in the Gibbons--Hawking form
\be
\label{GH_form}
g=Vh_{flat}+V^{-1}(dT+A)^2,\quad\mbox{where}\quad T\equiv x,
\ee
where $h_{flat}$ is a flat metric on $\R^{1, 2}$, and
$V$ and $A$ are respectively a function, and a one--form
on $\R^{1,2} $ which satisfy the Abelian monopole equation
\[
dV=*dA,
\]
where $*$ is the Hodge endomorphism of $h_{flat}$.
The function $V$ can be read--off directly from
$g$, and is given by
\[
V=g(K, K)^{-1}=\frac{ty-z^2}{2a^3t^3}.
\]
To construct the flat coordinates for the metric
\[
h_{flat}=V^{-1}(g-V K\otimes K)
\]
we first compute the six generators of its group of isometries. We then select a three-dimensional abelian subalgebra $(X_1, X_2, X_3)$ generated by translations. The corresponding one forms
$h_{flat}(X_i, .)$ are exact differentials of the flat coordinates $(Y, Z, T)$ where
\be
\label{old_coord}
y=\frac{X^2-Y^2-Z^2}{2(X+Y)^{3/2}}, \quad z=\frac{Z}{(X+Y)^{1/2}}, \quad t=(X+Y)^{1/2}.
\ee
Now
\be
\label{V_for_3}
h_{flat}=a^6 (dX^2-dY^2-dZ^2), \quad
V=\frac{X^2-Y^2-3Z^2}{4a^3(X+Y)^{5/2}}
\ee
and $V$ satisfies the wave equation on $\R^{1,2}$.
\begin{itemize}
\item
Instead of using the patching relation (\ref{patch_2}) we could have started with
\[
\tilde{\tau}=\lambda\tau+b\lambda^{-1}Q+a\lambda^{-2}Q^2,
\]
which allows the limit $a\rightarrow 0$ corresponding to the patching
(\ref{patch_flat}), and resulting in
a flat conformal structure. The corresponding
metric and the conformal factor arise from
$g$ and $\Omega$ given by (\ref{g_2}) and (\ref{conf_f}) when one makes a replacement
\[
z\longrightarrow z+\frac{b}{2a}.
\]
Therefore, if $a\neq 0$ then $b$ can be set to zero by
translating $z$.
\item The metric (\ref{g_2})
admits a second Killing vector
\[
K_2=5x\p_x+2y\p_y-z\p_z-4t\p_t
\]
which is not tri-holomorphic.
It Lie derives
$\Sigma^{01}$, but rotates $\Sigma^{00}$ and
$\Sigma^{11}$. Thus the space of orbits of $K_2$ in $M_\C$ admits a
Toda
Einstein--Weyl structure \cite{BF,Ward1}.
\item
There exists a  combination of self--dual  two forms which
is degenerate when the harmonic function $V$ in (\ref{g_2})
vanishes, which is the surface $ty-z^2=0$ in $M_\C$. All three
forms vanish on the line $t=y=z=0$ in $M_\C$. The normal bundle of the
twistor curves corresponding to this line jumps from
$\OO(1)\oplus\OO(1)$ to $\OO(3)\oplus\OO(-1)$.
\end{itemize}
\subsection{Ricci--flat metric
with $k=4$}
\label{k=4ex}
Parametrise the sections of ${\OO(4)}\rightarrow\CP^1$
by
\be
\label{parameter4}
Q=t+\lambda z+\lambda^2y+\lambda^3 w+\lambda^4 x.
\ee
In this case the splitting (\ref{split_2}) is possible if
\be
\label{cone}
\phi\equiv tw+zy=0.
\ee
Computing the resultant (\ref{resultant_g})
leads to the conformal structure
\[
g=\Omega^2\Big(t^6 dx^2-t^3z(ty+z^2)dxdz+2t^4(ty-z^2)dxdy+t^2(2t^2y^2-5tyz^2+z^4)dxdt-2tyz^2(ty-z^2)dz^2
\]
\be
\label{g_4}
-tz(t^2y^2-z^4)dzdy
-yz(t^2y^2-6tyz^2+z^4)dzdt+t^2(ty-z^2)^2dy^2
\ee
\[
+2t^2y^2(ty+z^2)dydt+y^2(t^2y^2+4tyz^2-z^4)dt^2\Big).
\]
The conformal factor making this metric Ricci flat is
\[
\Omega^2=\frac{2a^4}{t^2z(3ty-z^2)}.
\]
The triholomorphic Killing vector $\p/\p x$  Lie derives the
covariantly constant basis
of self-dual two--forms
\begin{eqnarray}
\label{sigma_k}
\Sigma^{00}&=&2 dx\wedge (tdy+ydt+zdz),\nonumber\\
\Sigma^{01}&=&tdx\wedge dz+zdx\wedge dt-\frac{yt-z^2}{t}dz\wedge dy-
\frac{y(yt+z^2)}{t^2}dz\wedge dt-\frac{2yz}{t}dy\wedge dt,\nonumber\\
\Sigma^{11}&=&dt \wedge\Big(4\frac{yz}{t}dz-2{t}dx-
2\frac{yt-z^2}{t}dy\Big).
\end{eqnarray}
The coordinate transformation (\ref{old_coord}) brings
the metric $g$ to Gibbons--Hawking form (\ref{GH_form}),
where
\be
\label{V4}
V=\frac{z(3ty-z^2)}{2a^4t^4}=
\frac{Z(3X^2-3Y^2-5Z^2)}{4a^4 (X+Y)^{7/2}},
\ee
and
\[
h_{flat}=a^8(dX^2-dY^2-dZ^2).
\]
The metric $g$ admits a homothety, as well as a
second non-triholomorphic
Killing vector $K_2=t\p_t-y\p_y-3x\p_x$.
\subsection{Gibbons--Hawking potential for general k} There exists a  map from  an affine bundle over $\OO(k)$
with holomorphic charts $(Q, \tau)$ and
$(\widetilde{Q}, \tilde{\tau})$
to an affine bundle over $\OO(2)$ with charts $(q, p)$ and $(\tilde{q}, \tilde{p})$
such that the latter admits a four--parameter family of section only if the patching for the former satisfies some additional conditions. The explicit transformation is given by
\[
q=\lambda^k\tau+Q^2, \quad p=\sum_{n=1}^{\infty}
\frac{(2n)!}{(1-2n)(n!)^2 4^n}\lambda^{(n-1)k}\tau^n(\lambda^k \tau+Q^2)^{1/2-n}
\]
on  $U$, and
\[ \tilde{q}=\tilde{\tau}, \quad \tilde{p}=\widetilde{Q}.
\]
on $\widetilde{U}$. This map is well defined only if some
sections are removed from the $\OO(k)$ twistor space. This corresponds to removing the region from $M_\C$ corresponding to the `big jump'. In the case of (\ref{patch_2})
we find
\[
\tilde{q}=\lambda^{-2}q, \quad \tilde{p}=p+
s(q, \lambda), \quad\mbox{where}\quad
s=\lambda^{-k}\sqrt{q}.
\]
The element of $H^1(\OO(2), \OO(-2))$ corresponding to the Gibbons--Hawking
function
is $\p s /\p q$. Parametrising the sections of $\OO(2)$ by
\[
q=\lambda^2 (X-Y)+ 2\lambda Z+ (X+Y)
\]
and taking the contour enclosing $\lambda=0$ in the twistor integral formula, leads to
\[
V(X, Y, Z)=\frac{1}{2(k-1)!}\frac{\p^{k-1}}{\p \lambda^{k-1}}
(\lambda^2 (X-Y)+2\lambda Z+(X+Y))^{-1/2}|_{\lambda=0},
\]
which for $k=3$ and $k=4$ agrees with (\ref{V_for_3}) and (\ref{V4}).
\subsection{$k=2$ and $k=1$} For completeness we shall analyse the remaining cases
$k=2$ and $k=1$ with the constant $a$ set to $1$. If $k=2$, both sides of
(\ref{split_one}) are homogeneous of degree $0$, and thus are equal to
some $x_{-1}$, so that
\[
\tilde{\tau}=x_{-1}+2\tilde{\lambda} x_0x_1+\tilde{\lambda}^2 x_0^2,
\quad \widetilde{Q}=x_2+\tilde{\lambda}x_1+\tilde{\lambda}^2 x_0,
\]
where in $\tilde{\tau}$ we have absorbed $(2x_0x_2+{x_1}^2)$ into $x_{-1}$.
Now compute the resultant (\ref{resultant_g}), and set
\[
x_{-1}=-2iT+\frac{1}{2}(X^2+Y^2)-Z^2, \quad x_0=\frac{1}{\sqrt{2}}(X+iY), \quad x_1=\sqrt{2}iZ,
\quad x_2=\frac{1}{\sqrt{2}}(X-iY).
\]
This yields the folded hyper-K\"ahler
metric (\ref{folded_gg}).
\vskip3pt
Finally if $k=1$ then both sides of (\ref{split_one}) are homogeneus of degree $1$, and thus
give rise to two coordinates on $M_\C$. The resulting metric is flat. In this case
the twistor space $\mathcal{Z}$ fibres over $\OO(1)$, and all metrics (corresponding to arbitrary patching)
fall into the classification of \cite{DW}.

\section{From Sparling--Tod to Eguchi--Hanson}
\label{EHST}
The holomorphic Sparling--Tod metric
\cite{Sparling_Tod,Tod1}
\be
\label{st_metric}
g=4dudv-4dxdy-8\rho\triangle^{-3}(udv-xdy)^2, \qquad
\triangle\equiv uv-xy
\ee
is ASD,  Ricci--flat, and of Petrov-Penrose type
$D$. In \cite{TN} a twistor--theoretic argument was
used to show that there exists
a Riemannian real section of (\ref{st_metric}) which is equivalent to the Eguchi--Hanson gravitational instanton. The coordinate transformation below
makes this explicit, by putting
(\ref{st_metric}) in the ALE $A_2$ Gibbons--Hawking form.

We find that the four--dimensional isometry group
of (\ref{st_metric}) contains $SL(2, \C)$ which acts tri--holomorphically. Let us consider a pencil of tri--holomorphic Killing vectors given by
\[
K=\frac{b}{2}(v\p_y+x\p_u)-\frac{c}{2}(y\p_v+u\p_x),
\]
where $(b, c)$ are constants not both zero.
A parallel basis of ${\Lambda^2}_+$ is Lie--derived along $K$, and the corresponding moment maps are
the flat coordinates on $\R^3$ in the Gibbons--Hawking form.
They are given by
\[
Z=i(bxv+cyu), \quad
X+iY=\sqrt{2}\rho(bx^2+cu^2)\triangle^{-2}+
\frac{\sqrt{2}}{2}(bx^2+cu^2), \quad
X-iY=\sqrt{2}(cy^2+bv^2).
\]
The metric (\ref{st_metric}) takes the form
\[
g=V(dX^2+dY^2+dZ^2)+V^{-1}(dT+A)^2,
\]
where $V$ is the harmonic function on $\R^3$ given by
\[
V=\frac{\triangle^{3}}{2\rho Z^2+bc\triangle^{4}}
\]
and
\be
\label{triangle_1}
\triangle^2=\frac{2bc\rho-
R^2 +\epsilon\sqrt{(R^2-2bc\rho)^2+8bc\rho Z^2}}{2bc}, \quad R^2\equiv X^2+Y^2+Z^2
\ee
where $\epsilon=\pm 1$.
With the help of some algebra this simplifies to
\[
V=\frac{1}{2\sqrt{-bc}}\;\Big(\frac{1}{|{\bf R}+{\bf a}|}-\epsilon \frac{1}{|{\bf R}-{\bf a}|}\Big),
\]
where ${\bf R}=(X, Y, Z), {\bf a}=(0, 0, \sqrt{-2bc\rho})$.
If $b, c$ are real and such that
$bc<0$ then $V$ with $\epsilon=-1$ gives the positive--definite Eguchi--Hanson gravitational instanton.
If $\epsilon=1$ then $g$ is still positive--definite, but not complete. It is an example
of a folded hyper--K\"ahler metric \cite{H2}. The points on the hypersurface $V=0$ in $M$ correspond to lines in $\mathcal{Z}$ where the normal bundle jumps
to $\OO\oplus \OO(2)$. Moreover the limit when
$b$ or $c$ tends to zero gives a dyon.
\section{Multi jumps}
\label{sec_casca}
The metrics resulting from Theorem \ref{main_theo} admit a tri--holomorphic Killing vector,
and thus can locally be put in the holomorphic Gibbons-Hawking form \cite{GH}, which depends on a
solution  $V$  to a holomorphic Laplace equation on
$\C^3$. The hypersurface corresponding to $V=0$ is singular, and can be characterised by
the jumping phenomenon.
\begin{prop}
Let $(M_\C, g)$ be a Gibbons--Hawking metric
\[
g=V(dX^2-dY^2-dZ^2)+V^{-1}(d T+A)^2, \quad\mbox{where}\quad
dV=\ast_3 dA
\]
and let $S=\{p\in M_\C, V(p)=0\}$. The points of $S$ correspond to
rational curves in $\mathcal{Z}$ with normal bundle $\OO(2)\oplus\OO$.
\end{prop}
{\bf Proof.}
The twistor space of a Gibbons-Hawking manifold
is an affine line bundle over the total space of $\OO(2)$
with transition functions
\[
\widetilde{\tau}=\tau+f(Q,\lambda), \quad
\widetilde{Q}=\lambda^{-2}Q,
\]
where $f\in H^{1}(\OO(2), \OO)$.
Restricting the cohomology class $f$ to a section
of $\OO(2)$
\be
\label{ofunction}
Q=\lambda^2 (X-Y)+2\lambda Z + (X+Y)
\ee
gives rise to the harmonic function $V$ by
\[
V(p)=\frac{1}{2\pi i}
\oint_{\Gamma\subset L_p}\frac{\p f}{\p Q} d\lambda.
\]
The normal bundle to $L_p$ is the restriction to
(\ref{ofunction}) of
\[\left(\begin{array}{cc}
           1 & \frac{\partial f}{\partial Q} \\
          0 & \lambda^{-2}\\
\end{array}\right)
\]
and then expanding
\[
f(Q,\lambda)=\sum_{-\infty}^{\infty}a_k\lambda^k, \quad\mbox{with}\quad
a_k=a_k(X, Y, Z).
\]
Split the sum into two:
\[\tilde{h}=-\sum_{-\infty}^{-1}a_k\lambda^k,\quad
h=\sum_0^{\infty}a_k\lambda^k,\quad\mbox{so that}\quad
f=h-\tilde{h}
\]and there is freedom to add a function
of $(X, Y, Z)$ to both of $h, \tilde{h}$.
Note that, after restricting,
\[\frac{\partial f}{\partial Z} =\frac{\partial f}{\partial Q} \frac{\partial Q}{\partial Z}=2\lambda\frac{\partial f}{\partial Q}\]
so that the transition matrix for the normal bundle is
\[F:=\left(\begin{array}{cc}
           1 & \frac{1}{2\lambda}\frac{\partial f}{\partial Z} \\
          0 & \lambda^{-2}\\
\end{array}\right).
\]
This is equivalent to
\[F\rightarrow \left(\begin{array}{cc}
           1 & \tilde{p} \\
          0 & 1\\
\end{array}\right)F\left(\begin{array}{cc}
           1 & p \\
         0 & 1\\
\end{array}\right)=\left(\begin{array}{cc}
           1 &  \frac{1}{2\lambda}\frac{\partial f}{\partial Z}+p+\tilde{p}\lambda^{-2}\\
          0 & \lambda^{-2}\\
\end{array}\right)\]
where we choose $\tilde{p},p$ to remove from $\frac{1}{2\lambda}\frac{\partial f}{\partial Z}$ all non-negative powers of $\lambda$ and all negative powers less than or equal to $-2$. All that remains is $\frac{1}{2\lambda}\frac{\partial a_0}{\partial Z}$, and
$\frac{\partial a_0}{\partial Z}$ is equal to a multiple of $V$.
Where $V\neq 0$, we know that this is the transition matrix for $\OO(1)\oplus\OO(1)$ but clearly where $V=0$ it is the transition matrix for $\OO(2)\oplus \OO$.
\koniec

We conclude that the metrics arising from Theorem \ref{main_theo} have at least two jumps:
the `small jump' from $\OO(1)\oplus\OO(1)$ to $\OO(2)\oplus \OO$ on a surface
$S$ corresponding to the zero set of $g(\p_x, \p_x)^{-1}$, and a big jump
to $\OO(k)\oplus\OO(2-k)$ on a curve $\gamma$. The argument below demonstrates that
many intermediate jumps can arise.

\vskip5pt
Consider the twistor space of Theorem \ref{main_theo}, with the moduli space
of rational curves $M_\C$ given by the Zariski cone (\ref{cones}). The constraints
defining $M_\C$ take the form
\[
x_0x_n+x_1x_{n-1}+\ldots =0, \quad\mbox{for}\;\; 3\leq n\leq k-1.
\]
For even $n=2m$ the last term in the constraint is $x^2_m/2$ and there will be a constraint like this for
\[
1\leq m\leq (k-1)/2 \quad\mbox{for odd}\; k \quad\mbox{or} \quad
1\leq m\leq k/2-1 \quad \mbox{for even}\; k.
\]
We shall be interested in solutions of the constraints for which all but one $x_n$
are zero (for $n<k$) and the constraints will not be satisfied for $n$ below a threshold. Thus we have a range of allowed $n$, namely $(1+k)/2\leq n\leq k-1$ for odd $n$ or $k/2\leq n\leq k-1$ for $n$ even.
For $n$ in these ranges the constraints are
satisfied with $x_n\neq 0$ and $x_i=0$ for all other $i$ in the range $0\leq i\leq k-1$. With any one of these solutions of the constraints, multiply $F$ on the right with
\be
\label{H_move}
H^{-1}=\left(\begin{array}{cc}
           1 & -2ax_k \\
          0 & 1\\
\end{array}\right)
\ee
to remove $x_k$-term from $Q$, leaving
\[
F=\left(\begin{array}{cc}
           \lambda^{k-2} & 2a\sum_{i} x_i\lambda^{i-2} \\
          0 & \lambda^{-k}\\
\end{array}\right).
\]
Now consider the product
\[\widetilde{H}\left(\begin{array}{cc}
           \lambda^{i-2} & 0 \\
           0 & \lambda^{-i} \\
\end{array}\right)H^{-1}=
\left(\begin{array}{cc}
           \alpha & 0 \\
          \lambda^{2-i-k} & -\alpha^{-1}\\
\end{array}\right)\left(\begin{array}{cc}
           \lambda^{i-2} & 0 \\
           0 & \lambda^{-i} \\
\end{array}\right)\left(\begin{array}{cc}
           \alpha^{-1}\lambda^{k-i} & 1 \\
          -1 & 0\\
\end{array}\right),
\]
when $\widetilde{H}$ and $H$ respectively have only negative or only positive powers of $\lambda$ and this product is $F$ given the
choice $\alpha=2ax_i$. Thus the normal bundle has jumped to $\OO(2-i)\oplus\OO(i)$ and there is an example like this for each $i$ in the allowed range.

We also always have the case $x_0\neq 0$, other $x_n$ zero when
\[\left(\begin{array}{cc}
          \alpha & 0 \\
         \lambda^{2-k} & -\alpha^{-1}\\
\end{array}\right)\left(\begin{array}{cc}
           \lambda^{-2} & 0 \\
           0 & 1 \\
\end{array}\right)\left(\begin{array}{cc}
           \alpha^{-1}\lambda^{k} & 1 \\
          -1 & 0\\
\end{array}\right)=\left(\begin{array}{cc}
           \lambda^{k-2} & \alpha\lambda^{-2} \\
          0 & \lambda^{-k}\\
\end{array}\right)=F
\]
with $\alpha=2ax_0$, so the jump to $\OO\oplus\OO(2)$ is always present.

We will see that all jumps are present if $k=4$. There is enough here to prove this also for
$k=5$ and $k=6$ but there is a gap at $k=7$: the above constructions do not give an example of a curve with normal bundle  $\OO(-1)\oplus\OO(3)$ but everything else up to $\OO(-5)\oplus\OO(7)$ occurs.
\subsection{Jump cascade with $k=4$.}
\label{jump_cascade}
Restrict the transition function (\ref{patch_2}) with $k=4$ to the line
(\ref{parameter4}).
The normal bundle is $\OO(1)\oplus\OO(1)$ away from $V=0$, where $V$ is given by (\ref{V4}),
and  jumps to $\OO(-2)\oplus\OO(4)$ at $y=z=w=t=0$. We have to look at the zero-set of $V$.
Multiply $F$ on the right with (\ref{H_move}) leaving
\[
F=\left(\begin{array}{cc}
           \lambda^2 & 2a(w\lambda+y+\frac{z}{\lambda}+\frac{t}{\lambda^2}) \\
          0 & \lambda^{-4}\\
\end{array}\right).
\]
There are six loci to investigate
all of which have $V=0$.
\begin{enumerate}
\item[$S_1$.] $w=y=z=t=0$ when we know it jumps to $\OO(-2)\oplus\OO(4)$.
\item[$S_2$.] $w=y=z=0$ but $t\neq 0$. Note that
\[
H^{-1}=\left(\begin{array}{cc}
           1/\beta & 0 \\
          \lambda^4 & \beta\\
\end{array}\right),\quad
\widetilde{H}=\left(\begin{array}{cc}
           0 & 1 \\
          -1 & \frac{1}{\beta\lambda^2}\\
\end{array}\right)
\]
give
\[\widetilde{H}\left(\begin{array}{cc}
           1 & 0 \\
          0 & \lambda^{-2}\\
\end{array}\right)H^{-1}=\left(\begin{array}{cc}
           \lambda^2 & \beta\lambda^{-2} \\
          0 & \lambda^{-4}\\
\end{array}\right)
\]
which with $\beta=2at$ shows $S_2$ is $\OO\oplus \OO(2)$.
\item[$S_3$.] $y=z=t=0$ but $w\neq 0$.
Take
\[
H^{-1}=\left(\begin{array}{cc}
           \lambda/\beta & 1 \\
          -1 & 0\\
\end{array}\right),\quad
\widetilde{H}=\left(\begin{array}{cc}
           \beta & 0 \\
          \lambda^{-5} & \frac{1}{\beta}\\
\end{array}\right)
\]
which give
\[\widetilde{H}\left(\begin{array}{cc}
           \lambda & 0 \\
          0 & \lambda^{-3}\\
\end{array}\right)H^{-1}=\left(\begin{array}{cc}
           \lambda^2 & \beta\lambda \\
          0 & \lambda^{-4}\\
\end{array}\right)
\]
which with $\beta=2bw$ shows $S_3$ is $\OO(-1)\oplus\OO(3)$.
\item[$S_4$.] $z=t=0$ but $y\neq 0$ (any $w$).
Take
\[
H^{-1}=\left(\begin{array}{cc}
           \lambda^2 & \alpha+\beta\lambda \\
          \frac{\beta\lambda}{\alpha^2}-\frac{1}{\alpha} & \frac{\beta^2}{\alpha^2}\\
\end{array}\right),\quad \widetilde{H}=\left(\begin{array}{cc}
           1 & 0 \\
          -\frac{\beta}{\alpha^2\lambda^3}+\frac{1}{\alpha\lambda^4} & 1\\
\end{array}\right)
\]
which  give
\[
\widetilde{H}\left(\begin{array}{cc}
           1 & 0 \\
          0 & \lambda^{-2}\\
\end{array}\right)H^{-1}=\left(\begin{array}{cc}
           \lambda^2 & \alpha+\beta\lambda \\
          0 & \lambda^{-4}\\
\end{array}\right)
\]
so with $\alpha=2ay,\beta=2az$ this shows that $S_4$ is $\OO\oplus\OO(2)$.
\item[$S_5$.] $z=w=0$ with $yt\neq 0$. Consider
\[
H^{-1}=\left(\begin{array}{cc}
            -\frac{\lambda^2}{\alpha}+\frac{\beta}{\alpha^2} & -1 \\
          1 & 0\\
\end{array}\right),\quad \widetilde{H}=\left(\begin{array}{cc}
           -\alpha -\frac{\beta}{\lambda^2} & \frac{\beta^2}{\alpha^2} \\
          -\frac{1}{\lambda^4} & -\frac{1}{\alpha}+\frac{\beta}{\alpha^2\lambda^2}\\
\end{array}\right)
\]
so that
\[\widetilde{H}\left(\begin{array}{cc}
           1 & 0 \\
          0 & \lambda^{-2}\\
\end{array}\right)H^{-1}=\left(\begin{array}{cc}
           \lambda^2 & \alpha+\frac{\beta}{\lambda^2} \\
          0 & \lambda^{-4}\\
\end{array}\right).
\]
With $\alpha=2ay,\beta=2at$ and $yt\neq 0$ this shows that $S_5$ is $\OO\oplus\OO(2)$.
\item[$S_6$.] $z^2=3ty$ with $yzt\neq 0$.
Introduce
\[\chi:=\frac{t}{\lambda^2}+\frac{z}{\lambda}+y+w\lambda\]
then set $t=\beta,z=3\alpha t$. We have  $3yt-z^2=0=wt+yz$ so that $y=3\alpha^2\beta$ and $w=-9\alpha^3\beta$ whence
\[\chi=\frac{\beta}{\lambda^2}(1+3\alpha\lambda+3\alpha^2\lambda^2-9\alpha^3\lambda^3)\]
and this is the top-right entry in $F$. Consider
 \[
H^{-1}=\left(\begin{array}{cc}
           1-3\alpha\lambda+6\alpha^2\lambda^2 & 9\alpha^4\beta(5-6\alpha\lambda) \\
          \frac{\lambda^4}{\beta} & \frac{\lambda^2f}{\beta}\\
\end{array}\right),\quad \widetilde{H}\left(\begin{array}{cc}
           0 & \beta \\
          -\frac{1}{\beta} & \frac{1}{\lambda^2}-3\frac{\alpha}{\lambda}+6\alpha^2\\
\end{array}\right).
\]
Then
\[\widetilde{H}\left(\begin{array}{cc}
           1 & 0 \\
          0 & \lambda^{-2}\\
\end{array}\right)H^{-1}=\left(\begin{array}{cc}
           \lambda^2 & \chi \\
          0 & \lambda^{-4}\\
\end{array}\right)
\]
so $S_6$ is also $\OO\oplus\OO(2)$.
\end{enumerate}
\medskip
We conclude that there are curves with normal bundle $(1,1),(0,2),(-1,3)$ and $(-2,4)$ i.e. all possibilities up to the maximum jump occur.

\section{Generalised Legendre transform and self--dual two--forms}
\label{sec_leg}
There is another route directly from the cohomology class defining the affine line bundle
$\mathcal{Z}\longrightarrow {\mathcal O}(k)$ to the ASD Ricci--flat metric directly
without the need to use resultants as in (\ref{resultant_g}). This follows
Theorem 4.4 in \cite{DM}, and gives a version of the generalised
Legendre transform \cite{IR96,B00}.

Affine line bundles over ${\mathcal O}(k)$ are classified by elements $[f]$ of $H^1(\OO(k),
\OO(2-k))$ as
\be
\label{affine_line_b}
\tilde{\tau}=\tau+f(Q, \lambda), \quad \widetilde{Q}=\lambda^{-k}Q.
\ee
Any such cohomology class gives rise to
$k-3$ constraints
\be
\label{constrain4}
\phi_{A_1\dots
A_{k-4}}:=\frac{1}{2\pi i}\oint_{\Gamma}\pi_{A_1}\dots \pi_{A_{k-4}}f(Q,
\pi_{A})\pi\cdot d\pi=0
\ee
which trace out a holomorphic four--manifold $M_\C$ in a $(k+1)$--dimensional space of holomorphic
sections of $\OO(k)\rightarrow \CP^1$. The ASD Ricci--flat metric on $M_\C$ is determined by a  basis
of self--dual two forms $\{\Sigma^{00}, \Sigma^{01}, \Sigma^{11}\}$ which are pull-backs
from $\C^{k+1}$ to $M$ of two--forms\footnote{This formula
corrects (4.27) in \cite{DM}.}
\begin{eqnarray}
\label{SDtwoform}
\Sigma^{AB}&=&
\frac{1}{8}{\psi^{AB}}_{B_1\dots B_{k-3}C_1\dots C_{k-3}} d {x_{PQR}}^{B_1\dots B_{k-3}}\wedge
 d {x^{PQRC_1\dots C_{k-3}}}+\\
&&\frac{3}{2}
\psi_{B_1\dots B_{k-2}C_1\dots C_{k-2}}
{dx_P}^{B_1\dots B_{k-2}(A}\wedge dx^{C_1)C_2\dots C_{k-2}BP},\nonumber
\end{eqnarray}
where
\be
\label{the_field}
\psi_{A_1\dots A_{2k-4}}=\frac{1}{2\pi i}\oint_{\Gamma}
\pi_{A_{1}}\dots\pi_{A_{2k-4}}
\frac{\p f}{\p Q}\pi\cdot d\pi
\ee
is a zero--rest--mass field determined by $[f]$, and $Q$ in this formula is regarded as the coordinate on the fibres of $\OO(k)\rightarrow\CP^1$
which is homogeneous of degree $k$. 

If $k=3$ then there are no constraints
to be imposed,  and $\psi_{AB}$ is a self--dual Maxwell field originally constructed
in \cite{W78}.
\subsection{Example with $k=4$}
The manifold
$M_\C$ is a surface $\phi=0$ given by (\ref{constrain4})  in the five--dimensional space
${\mathcal N}$ of holomorphic sections of the fibration $\OO(4)\rightarrow \CP^1$.
The function $\phi$ satisfies the overdetermined system of linear PDEs
${\p^{A}}_{BCD}\psi_{EFGA}=0$ where $\psi_{ABCD}$ is given by
(\ref{the_field}). Explicitly
\begin{eqnarray*}
&&\phi_{yt}-\phi_{zz}=0, \quad \phi_{tw}-\phi_{yz}=0, \quad
\phi_{tx}-\phi_{wz}=0,\\
&& \phi_{wz}-\phi_{yy}=0,\quad \phi_{xz}-\phi_{wy}=0, \quad \phi_{xy}-\phi_{ww}=0.
\end{eqnarray*}

 Consider the cohomology class represented by
$f=Q^2\lambda^{-k}$, and take $k=4$.
Comparing $Q=x^{ABCD}\pi_A\pi_B\pi_C\pi_D$
with (\ref{parameter4}) gives
\[
t=x^{1111},\quad z=4x^{1110}, \quad y=6x^{1100},
\quad w=4x^{1000},  \quad
x=x^{0000}.
\]
Evaluating the residue at the pole $\lambda=0$ in
(\ref{constrain4}) yields the constraint
\[
\phi=tw+zy=0
\]
in agreement with (\ref{cone}). The spin--2 field
(\ref{the_field}) is
\[
\psi_{0000}=0, \quad
\psi_{0001}=t, \quad \psi_{0011}=z,\quad
\psi_{0111}=y,\quad \psi_{1111}=w
\]
which gives the self-dual two forms
\begin{eqnarray*}
\Sigma^{00}&=&z(2dx\wedge dz+\frac{1}{2}dw\wedge dy)+
y(2dx\wedge dt-\frac{1}{2}dz\wedge dw)
+w(\frac{1}{2}dw\wedge dt)+t(2dx\wedge dy),\\
\Sigma^{01}&=&z(dx\wedge dt-dz\wedge dw)+
y(dw\wedge dt+\frac{1}{2}dy\wedge dz)
+w(\frac{1}{2}dy\wedge dt)+t(dx\wedge dz+\frac{1}{2}dw\wedge dy),\\
\Sigma^{11}&=&z(2dw\wedge dt+\frac{1}{2}dy\wedge dz)+
y(2dy\wedge dt)
+w(\frac{3}{2}dz\wedge dt)+t(2dx\wedge dt -\frac{1}{2}dz\wedge dw).
\end{eqnarray*}
The pull-back of these two--forms to the cone
(\ref{cone}) agrees with expressions (\ref{sigma_k}).
\vskip5pt
The jump cascade discussed in Section \ref{jump_cascade} can be now understood in the framework of the generalised Legendre transform presented in 
\cite{moraru} and \cite{DT18}.
 Using the Kodaira isomorphism
\[
{T^*}_p{\mathcal N}\cong  H^0(L_p, \OO(4))=\mbox{Sym}^4(\C^2)
\]
we can identify the gradient $d\phi$ with a binary quartic
\begin{eqnarray*}
d\phi \rightarrow  {\mathcal{Q}}(d\phi)&=&\alpha s^4+4\beta s^3+6\gamma
s^2+4\delta s+\epsilon\\
&=&\phi_x s^4-4\phi_w s^3+6\phi_y s^2-4\phi_z s+\phi_t.
\end{eqnarray*}
Binary quartics admit two classical invariants
\be
\label{inv_ij}
{\mathcal I}= \alpha\epsilon-4\beta\delta+3\gamma^2 , \quad\mbox{and}\quad
{\mathcal J}=\mbox{det} \left(\begin{array}{ccc}
          \alpha & \beta &\gamma \\
          \beta &\gamma &\delta  \\
           \gamma & \delta &\epsilon
\end{array}\right).
\ee
If $\phi$ is given by (\ref{cone}) then
\[
\mathcal{I}=3z^2-4ty, \quad
\mathcal{J}=z^3-2tzw+t^2w.
\]
The points in $M_\C$ where $d\phi=0$
correspond to twistor curves
with normal bundle $\OO(-2)\oplus\OO(4)$. The points where
$d\phi\neq 0$, but
 ${\mathcal I}={\mathcal J}=0$ correspond to twistor curves
with normal bundle $\OO(-1)\oplus\OO(3)$. The points
where ${\mathcal I}\neq 0$ and ${\mathcal J}=0$ correspond to curves with  normal bundle $\OO\oplus\OO(2)$. Finally the generic points have ${\mathcal I}\neq 0,  {\mathcal J}\neq 0$. Such points
correspond to twistor curves with the normal bundle
$\OO(1)\oplus\OO(1)$.
\subsection{A Riemannian example}
The Riemannian reality conditions require $k=2n$ to be even. The real sections
satisfy
\[
\overline{Q(\lambda)}=(-1)^n\ov{\lambda}^{2n}
Q(-1/\ov{\lambda}),
\]
which in the case $k=4$ imples  that
\[
Q=t+\lambda z+\lambda^2y-\lambda^3\ov{z}+\lambda^4\ov{t}
\]
with $t,  z$ complex  and $y$ real. The surface (\ref{cone}) becomes $t\ov{z}+zy=0$ which is of codimension
$2$ in the space of real sections of $\OO(4)$. Thus the metric (\ref{g_4}) does not  admit a Riemannian slice.

To construct a Riemannian metric which admits a jump to $\OO(-2)\oplus\OO(4)$
consider a twistor space defined by the patching
relation\footnote{To make contact with (\ref{affine_line_b})
divide the expression (\ref{patch_r}) for $\tilde{\tau}$ by
$\lambda^2$, and set $f=s/\lambda^2$.}
\be
\label{patch_r}
\widetilde{Q}={\lambda^{-4}}Q,\quad
\tilde{\tau}=\lambda^{2}\tau+s(Q, \lambda), \quad\mbox{where}\quad
s=3Q^2(1-\lambda^{-6}).
\ee
The metric can be computed
as in (\ref{section_zari})
using the resultant (\ref{resultant_g}), and constructing a conformal factor
which makes the metric Ricci--flat.
We shall instead perform the Legendre
transform of \cite{IR96} which  leads directly to a K\"ahler potential for the metric. To make contact with the notation and formalism of
\cite{IR96}  define $G(Q, \lambda)$ by
\[
\frac{\p G}{\p Q}=\frac{s}{\lambda^2}, \quad\mbox{so that}\quad
G=\frac{Q^3}{\lambda^2}(1-\lambda^{-6}),
\]
and set
\begin{eqnarray*}
F&=&\frac{1}{2\pi i}\oint_{\Gamma\subset\CP^1} \frac{1}{\lambda^2}G(t+\lambda z+\lambda^2y-\lambda^3\ov{z}+
\lambda^4\ov{t}, \lambda)d\lambda\\
&=&6ytz+ 6y\ov{tz}+z^3+\ov{z}^3-3z\ov{t}^2-3\ov{z}t^2,
\end{eqnarray*}
where the contour $\Gamma$ encloses $\lambda=0$.
The real four--manifold $M$ is defined as the surface
\begin{eqnarray*}
\phi:=\frac{\p F}{\p y}&=&6(tz+\ov{tz})\\
&=&0
\end{eqnarray*}
in the space of real sections of $\OO(4)$.
Using the splitting method
in the proof of Theorem
 \ref{main_theo}, or equivalently
computing the ${\mathcal{I}}$ and ${\mathcal{J}}$ invariants
(\ref{inv_ij}) we find that the points in $M$ where
$t=z=0$ correspond to curves with normal bundle
$\OO(-2)\oplus\OO(4)$. This is a curve parametrised by $y$.

Now perform the Legendre transform
\[
u:=\frac{\p F}{\p z}=6yt+3z^2-3\ov{t}^2
\]
and eliminate the coordinates
$(z, \ov{z}, y)$ using $(t, \ov{t}, u, \ov{u})$ as holomorphic and
anti--holomorphic coordinates on $M$. The K\"ahler potential is
\begin{eqnarray*}
\Omega(t, \ov{t}, u, \ov{u})&=&F-uz-\ov{uz}\\
&=&-2(z^3+\ov{z}^3)\\
&=&-2i(t^3-\ov{t}^3)R^3, \quad\mbox{where}\quad R^2=-1
-\frac{\ov{u}t-u\ov{t}}{3(t^3-\ov{t}^3)}\in \R^+.
\end{eqnarray*}
The K\"ahler potential satisfies the 1st heavenly equation \cite{pleban}
\[
\Omega_{t\ov{t}}\Omega_{{u}{\ov{u}}}-\Omega_{t\ov{u}}\Omega_{{u}\ov{t}}=1
\]
and the resulting metric on $M$
\[
g=\Omega_{{u}\ov{u}}dud\ov{u}+
\Omega_{{u}\ov{t}}dud\ov{t}+\Omega_{{t}\ov{u}}dtd\ov{u}
+\Omega_{{t}\ov{t}}dtd\ov{t}
\]
is hyper--K\"ahler. The line of jumping points in $M$ has been blown down to a point $u=t=0$ by the Legendre transform.
\section{Schr\"odinger equation on folded hyper--K\"ahler manifolds}
\label{sec_shrod}
In this Section we shall demonstrate that the Schr\"odinger
equation on a canonical folded hyper--K\"ahler manifold
(corresponding to $k=2$ in Theorem \ref{main_theo})
\[
g=Z\left(dX^{2}+dY^{2}+dZ^{2}\right)+
\frac{1}{Z}\left(dT+\frac{1}{2}XdY-\frac{1}{2}YdX\right)^{2}
\]
admits normalisable solutions which extend to both sides of the fold $Z=0$ where the metric
degenerates.

The time-independent Schr\"odinger equation
\[
\frac{1}{\sqrt{|g|}}\partial_{a}\left(\sqrt{|g|}g^{ab}\partial_{b}\phi\right)=E\phi
\]
takes the form
\begin{equation}
\frac{1}{Z}(\frac{1}{4}(X^{2}+Y^{2})+Z^{2})\partial_{T}\partial_{T}\phi-\frac{X}{Z}\partial_{Y}\partial_{T}\phi+\frac{Y}{Z}\partial_{X}\partial_{T}\phi+\frac{1}{Z}\delta^{ij}\partial_{i}\partial_{j}\phi=E\phi.
\label{eq:SchOnFolded}
\end{equation}

We shall take the coordinate $T$ to be periodic, and consider
solutions of the form
\[
\phi(T, X, Y, Z)=e^{isT}\varphi(X, Y, Z)
\]
for $s$ a non-zero integer. The Schr\"odinger
equation (\ref{eq:SchOnFolded}) becomes
\[
-\frac{s^{2}}{Z}(\frac{1}{4}(X^{2}+Y^{2})+Z^{2})\varphi-\frac{isX}{Z}\partial_{Y}\varphi+\frac{isY}{Z}\partial_{X}\varphi+\frac{1}{Z}\delta^{ij}\partial_{i}\partial_{j}\varphi=E\varphi\,\,,
\]
which separates as $\varphi=G(X, Y)F(Z)$ into
\begin{equation}
\frac{d^{2}F}{dZ^{2}}-(s^{2}Z^{2}+EZ+\kappa)F=0\label{eq:Feqn}
\end{equation}
and
\begin{equation}
-\frac{1}{4}s^{2}(X^{2}+Y^{2})-\frac{isX}{G}\partial_{Y}G+\frac{isY}{G}\partial_{X}G+
\frac{1}{G}(\partial_{X}^{2}+\partial_{Y}^{2})G+\kappa=0.\label{Aeqn}
\end{equation}
If $s=0$ then the first equation becomes the Airy equations and one can show that no--normalisable
solutions exist on both sides of the fold.  The second equation describes
a free particle on a plane, and no bound states exist in this case either.

\vskip5pt
Let us therefore assume that $s\neq 0$, and consider
the equation for (\ref{eq:Feqn}) for $F(Z)$, which
has the form of the Schr\"odinger equation describing a displaced harmonic
oscillator. This is readily solved to give
\[
F(Z)=H_{\gamma}\left(\sqrt{s}\left(Z+\frac{E}{2s^{2}}\right)\right)\exp\left\{ -\frac{1}{2}s\left(Z+\frac{E}{2s^{2}}\right)^{2}\right\} ,
\]
where $H_{\gamma}(\xi)$ solves the Hermite equation
\[
\frac{d^{2}H}{d\xi^{2}}-2\xi\frac{dH}{d\xi}+2\gamma H=0,
\quad\mbox
{with}\quad
\gamma=\frac{1}{2s}\left(\frac{E^{2}}{4s^{2}}-(\kappa+s)\right).
\]
If $\gamma$ is a non-negative integer then $H_{\gamma}$ is a Hermite
polynomial and thus $F(Z)$ is clearly normalisable for $s>0$ (even
with the folded background's factor of $\sqrt{|g|}=Z$) due to the
exponential fall-off at large $Z$. If, however, $\gamma$ fails to
be a non-negative integer then $H_{\gamma}$ is more complicated,
being most readily expressed as a series expansion. In this case normalisability
is less clear, so let's restrict ourselves to the case where $\gamma$
is a non-negative integer.

\textcompwordmark{}

Let us now proceed to consider the $G(X ,Y)$ equation (\ref{Aeqn}).
This has the form of the Schr\"odinger equation describing motion in
a constant magnetic field. In the usual manner let us then define
the canonical (Hermitian) momenta
\[
\Pi_{X}=-i\partial_{X}+\frac{1}{2}sY\qquad\Pi_{Y}=-i\partial_{Y}-\frac{1}{2}sX
\]
and ladder operators
\[
a=\Pi_{X}+i\Pi_{Y}\qquad a^{\dagger}=\Pi_{X}-i\Pi_{Y}.
\]
The $G(X, Y)$ equation is then
\[
(a^{\dagger}a+s-\kappa)G=0
\]
and we can construct some solutions (choosing $\kappa=s$) by solving
$aG_{0}(X, Y)=0$, and then applying copies of $a^{\dagger}$ to $G_{0}$.
For example, one solution is
\[
G(X ,Y)\,\propto\,\exp\left\{ -\frac{1}{4}s(X^{2}+Y^{2})\right\} ,
\]
and thus we conclude that there do exist normalisable solutions. One
class of normalisable solutions is
\[
\phi=H_{\gamma}\left(\sqrt{s}\left(Z+\frac{E}{2s^{2}}\right)\right)\exp\left\{ -\frac{1}{2}s\left(Z+\frac{E}{2s^{2}}\right)^{2}\right\} \exp\left\{ -\frac{1}{4}s(X^{2}+Y^{2})\right\} \exp\left\{ isT\right\}
\]
with $s$ a positive non-zero integer and $E$ chosen such that
\[
\gamma=\frac{E^{2}}{8s^{3}}-1
\]
is a positive integer.

Another example of a metric which admits a three--parameter family of jumping lines,  and yet
there exists normalisable solutions to the Schr\"odiner equation is the Taub--NUT space with negative mass \cite{GM}.

\end{document}